\documentclass[11pt,twoside]{article}
\usepackage{amsmath}
\usepackage{amsthm}
\usepackage{amssymb}
\usepackage{array}
\usepackage{geometry}
\usepackage{graphicx}
\usepackage{enumerate}
\usepackage{lscape,graphicx}
\usepackage{fancyhdr}
\usepackage{makeidx}
\usepackage{amsfonts}
\usepackage{graphicx}
\usepackage{CJK}
\usepackage{amsmath}
\usepackage{latexsym}
\usepackage{amssymb}
\usepackage{mathrsfs}
\usepackage{amsthm}
\usepackage{array}
\usepackage{enumerate}
\usepackage{longtable}
\usepackage{multirow}
\makeindex
\usepackage{wasysym}
\usepackage{chapterbib}

\newtheorem{definition}{Definition}[section]
\newtheorem{thm}{Theorem}[section]

\normalsize
\begin{document}


\newpage
\setcounter{equation}{0}
\begin{center}
\vskip1cm
{\Large
\textbf{Maximization of Mathai's Entropy under the Constraints of Generalized Gini and Gini
mean difference indices and its Applications in Insurance} }
\vskip.5cm
%

\vspace{0.50cm}
\textbf{\large Rhea Davis and  Nicy Sebastian}\\
Department of Statistics,
St Thomas College, Thrissur, Kerala, India\\ {Email: {\tt{rheadavisc@gmail.com, nicycms@gmail.com}}}\\

\end{center}

\thispagestyle{empty}
\begin{center} {\bf }  \vskip 0.10truecm \noindent \\

\vskip.5cm\noindent \centerline{\bf Abstract}
\end{center}
\par
Statistical Physics, Diffusion Entropy Analysis and Information Theory commonly use Mathai's entropy which measures the randomness of probability laws, whereas welfare economics and the Social Sciences commonly use Gini index which measures the evenness of probability laws. Motivated by the principle of maximal entropy, we explore the maximization of Mathai's entropy subject to the conditions in the following scenarios: (i) the conditions of a density function and fixed mean; (ii) the conditions of a density function and fixed Generalized Gini index. We also maximizes the  Mathai's entropy subject to the constraints of a given Gini mean difference index and the conditions of a density function. The obtained maximum entropy distribution is fitted to the loss ratios (yearly data) for earthquake insurance in California from 1971 through 1994 and its performance with some one-parameter distributions are compared.
		 \\



\vskip.5cm{\section{Introduction}}
\vskip.3cm
The notion of entropy was first developed by physicists in the context of equilibrium thermodynamics. This concept was later extended to information theory and statistical mechanics. The most popular one is due to Shannon 
  and  he designated entropy as a measure of uncertainty or a measure of information, see Shannon (1948).
Various generalizations of Shannon entropy are available in the literature. R\'{e}nyi's entropy, Havrda-Charv\'{a}t entropy and Tsallis entropy are some of the important generalizations, see Mathai and Haubold (2007). 
 The maximum Tsallis distributions have encountered a huge success because of their remarkable agreement with experimental data. If the density function is replaced with escort density and if the expected value in this escort density is assumed to be fixed, then maximization of Tsallis entropy subject to these constraints will lead to Tsallis statistics, 
 see Tsallis (1988), 
Mathai and Haubold (2007). A notable addition to the collection of generalizations of Shannon entropy is the entropic form of order alpha or Mathai's entropy Mathai and Haubold (2006)
 ~ and Sebastian (2015)
. Mathai's entropy has been defined in cases where random variable X is a real scalar, complex scalar, real vector/matrix  and complex vector/matrix Mathai and Sebastian et al. (2021)
.

Approximation of distributions is a fundamental problem in statistical data analysis 
Tanak et al. (2015). The maximum entropy principle gives a general way of achieving this. The maximum entropy method for estimating the probability density function (pdf) was primarily suggested by
Jaynes (1957), and then it has been widely used in many research areas, mainly in economics.
Maximization of Shannon entropy and its generalizations subject to different sets of constrains are richly available in the literature, significantly in the field of economics and social sciences, in the context of estimating income distribution with regard to income inequality in society 
 see Kapur (1989), Kagan et al. (1973) and Tanak et al. (2017). In recent years, entropy maximization based on inequality measure constraints has been considered by several authors, see Liu et al. (2020), Rad et al. (2016), and Tanak et al. (2015).
According to this principle the distribution for inference should have the property that the density function maximizes entropy subject to certain constraints representing our incomplete information 
Tanak et al. (2017).

The paper is organized as follows. Section 2 contains some preliminaries and fundamental concepts of some of the important inequality measures, namely the generalized Gini index, available in the literature, which
will be used in the  sections that follow. In  Section 3, definitions, properties, and the maximization of Mathai's entropy under different sets of constraints available in the literature are presented.  In Section 4 Mathai's entropy is maximized subject to the conditions of a density function and fixed Generalized Gini Index.
In Section 5, we consider income distributions supported on the real line and derive
the maximum Mathai's entropy with mean and Gini mean difference index constraints.
In Section 6, the obtained maximum entropy distribution is fitted to the loss ratios (yearly data) for earthquake insurance in California from 1971 through 1994 and its performance with respect to some one-parameter distributions are compared.


\vskip.5cm{\section{Inequality Measures}}
\vskip.3cm
We need some method of knowing whether the distribution under consideration is becoming more or less unequal, see Lorenz (1905). 
 We wish to be able to say at what point a society is placed between the two extremes, equality, on the one hand, and the ownership of all wealth by one individual.  On the other hand inequality can be defined as the dispersion of the distribution of income or some other welfare indicator, see Tanak et. al (2017)
. The Lorenz curve introduced by Lorenz is an important tool for analyzing income inequality. Graphically, the Lorenz curve gives the proportion of total societal income flowing to the lowest earning proportion of the income earners.\\
	\begin{definition} [Lorenz curve]
		Let $X$ denote a random variable with cumulative distribution function (cdf) F supported in $(0,\infty)$ and mean $E(X)=\mu$. The Lorenz curve is defined as
		\begin{equation}\label{eq2}
			L_F(u) = \frac{1}{\mu}\int_{0}^{u}F^{-1}(x)dx, \, 0\leq u \leq 1,
		\end{equation}
		where $F^{-1}(x)=inf\{t:F(t)\geq x\}$.
	\end{definition}
	$ L_F(u) $ represents the income share of the total held by the lowest 100u\% of the population.\\

In discussions of appropriate fiscal policy, issues related to the distribution of income are important, see Dastrup (2007)
. At the core of the dialogue are the questions about who benefits from the changes in taxes and  transfer payments. In order to address these issues it is important to identify the model of the income distribution under consideration.\\
	Vilfredo Pareto first proposed a model of income distribution in 1895 in the form of a probability density function. This was found to be an accurate model in the case of the upper tail of the distribution, but did not perform well in the case of lower tail. Pareto's analysis sparked a debate on the effect of economic growth in income inequality. Gini disagreed with Pareto's ideas and introduced a measure of inequality known as the Gini index, which is an important measure of income inequality derived from the Lorenz curve.
	
	\begin{definition} [Gini Index]
		The Gini index, $G(F)$,  is defined as twice area between the considered Lorenz curve and the line of perfect equality $L_F(u)=u$.
		\begin{equation}\label{eq3}
			G(F) = 2\int_{0}^{1}(u-L_F(u))du = 1-2\int_{0}^{1}L_F(u)du.
		\end{equation}
	\end{definition}
	The Gini index takes values in the interval $[0,1]$. Greater the Gini index, greater is the income inequality in the society. The upper bound 1 indicates the situation where a single individual possesses all the income and all other members are completely impoverished.\\
	\hfill \\
	A single-parameter generalization of the Gini index is given by,
	\begin{equation}\label{eq4}
		G_\nu(F)=1-\int_{0}^{1}\nu(\nu-1)(1-u)^{\nu-2}L_F(u)du, \nu>1,
	\end{equation}
	where $ \nu $ is a parameter tuning the degree of 'aversion of inequality' so that higher weights are attached to smaller incomes as $ \nu $ increases. In case of $\nu=2$, we have Gini index. By using the definition of generalized Gini index (\ref{eq4}) and Lorenz curve (\ref{eq2}), we have
	\hfill \\
	\begin{equation}\label{eq5}
		G_\nu(F)=1-\frac{1}{\mu}\int_{0}^{\infty}{\bar{F}}^\nu(x)dx,
	\end{equation}
	where $\bar{F}(x)=1-F(x)$ is the survival function, for more details see Tanak et al. (2015).\\
Gini index is defined for probability distributions supported on the non-negative real line and hence, a measure of inequality for income distributions with real line support are needed. 

\begin{definition}[Gini mean difference (GMD) index]
Suppose a random variable $ X $ has a density function $ f(x), \, x \in \mathbb{R} $. Gini mean difference (GMD) index is given by,
\begin{equation} \label{gp1}
D(f) = \int_{-\infty}^{\infty} \int_{-\infty}^{\infty} |x-y| f(x)f(y) dxdy.
\end{equation}
\end{definition}
GMD is the expected absolute difference between two realizations of independently and identically distributed random variables, see Yitzhaki and Schechtman (2012) and Tanak et al. (2015). In terms of survival function $ \bar{F}(\cdot) $, $ D(f) $ can be represented as,
\begin{equation}\label{gp2}
D(f) = 2 \int_{-\infty}^{\infty} \bar{F}(t) F(t) dt.
\end{equation}
For more details see Gini (1912), Yitzhaki and Schechtman (2012) and Tanak et al. (2015).

\vskip.5cm{\section{ A Review of Mathai's Entropy}}
\vskip.3cm

In this section, definition and properties of Mathai's entropy are presented. The distributions obtained, including Mathai's pathway model, as a result of the
	maximization of Mathai's entropy under different sets of constraints, are also discussed.\\

\subsection{Definition and Properties}
	The generalized entropy of order $\alpha$ or Mathai's entropy  is defined as follows:
	\begin{definition}[Mathai's Entropy (Discrete case)] \label{d21}
		Consider a multinomial population $P = (p_1,p_2,...,p_k) , p_i\geq 0, i=1, ..., k, p_1+...+p_k=1$.  Then Mathai's entropy is defined as
		\begin{equation}\label{dismath}
			M_{k,\alpha}(P) = \frac{\sum_{i=1}^k p_i^{2-\alpha}-1}{\alpha-1}, \alpha\neq1, -\infty <\alpha<2.
		\end{equation}
	\end{definition}
	When $ \alpha $ decreases and is less than zero, then $ 1-\alpha $ increases and is greater than 1, and vice versa. Hence, $ 1-\alpha $ can be considered as the \textit{strength of information} in the distribution. Larger the value of $ 1- \alpha $ smaller the uncertainty and vice versa.

	\subsection{Properties}
	\begin{enumerate}
		\item{Non-negativity: $M_{k,\alpha}(P)\geq0$ with equality only when one $ p_i=1$ and the rest zeros.}
		\item{Expansibility or zero-indifferent: $M_{k+1,\alpha}(P,0)=M_{k,\alpha}(P)$. If an impossible event is incorporated into the scheme, that is, $ p_{k+1}=0$ it will not change the value of the entropy measure.}
		\item{Symmetry: $M_{k,\alpha}(P)$ is a symmetric function of $p_1,...,p_k$. Arbitrary permuations of $p_1,...,p_k$ will not alter the value of $ M_{k,\alpha}(P)$.}
		\item{Continuity: $M_{k,\alpha}(P)$ is a continuous function of $p_i>0, \,  i=1,...,k.$}
		\item{Monotonicity: $M_{k,\alpha}(\frac{1}{k},...,\frac{1}{k})$ is a monotonic increasing function of $k$.}
		\item{Inequality: $M_{k,\alpha}(p_1,...,p_k)\leq M_{k,\alpha}(\frac{1}{k},...,\frac{1}{k}).$}
		\item{Branching principle or recursivity:
			\begin{multline*}
				M_{k,\alpha}(p_1,...,p_k) = M_{k-1,\alpha}(p_1+p_2,p_3,...,p_k)+(p_1+p_2)^{2-\alpha}\\
				\times M_{2,\alpha}(\frac{p_1}{p_1+p_2},\frac{p_2}{p_1+p_2}).
			\end{multline*}
			This property indicates what happens to the measure if two of the mutually exclusive and totally exhaustive events are combined.}
		\item{Non-additivity: Consider independent multinomial populations $P= (p_1,...,p_n)$ and $Q = (q_1,..,q_m)$ such that $\sum_{i=1}^{n}\sum_{j=1}^{m}p_iq_j=1, \sum_{i=1}^{n}p_i=1,\sum_{j=1}^{m}q_j=1.$\\
			Then the joint density is of the form $(p_1q_1,...,p_1q_m,...,p_nq_1,...,p_nq_m)$. Let us denote the entropy measure in this joint distribution by $M_{nm,\alpha}(P,Q)$. Then
			\begin{multline*}
				M_{nm,\alpha}(P,Q) = M_{n,\alpha}(P)+M_{m,\alpha}(Q)+(\alpha-1)M_{n,\alpha}(P)M_{m,\alpha}(Q).
			\end{multline*}
			The third term on the right makes the measure non-addtive.}
		\item{Decomposibility: Consider the joint discrete distribution $ p_{ij}\geq 0,$\\$ \sum_{i=1}^n \sum_{j=1}^m p_{ij}=1$. Consider the marginal distribution $P_j = \sum_{i=1}^m p_{ij} > 0,  j=1,...,m$. Then we have
			\begin{multline*}
				M_{nm,\alpha}(p_{11},p_{12},...,p_{1m},p_{21},...,p_{2m},...,p_{n1},...,p_{nm})\\
				= M_{m,\alpha}(P_1,...,P_m)+\sum \limits_{j=1}^{m}P_j^{2-\alpha}M_{n,\alpha}(\frac{p_{1j}}{P_j},...,\frac{p_{nj}}{P_j}).
		\end{multline*}}
		\item{Functional equation: Consider$ M_{2,\alpha}(P)=M_{2,\alpha}(p,1-p).$ That is,
			\begin{center}
				$M_{2,\alpha}(p,1-p)=\frac{p^{2-\alpha}+(1-p)^{2-\alpha}-1}{\alpha-1}, \alpha\neq 1, \alpha < 2.$
			\end{center}
			Let $f_\alpha (p) = M_{2,\alpha}(p,1-p).$ Then $f_\alpha (p)$ satisfies the functional equation
			\begin{equation*}
				f_\alpha(x) + (1-x)^{2-\alpha}f_\alpha(\frac{y}{1-x})=f_\alpha(y)+(1-y)^{2-\alpha}f_\alpha(\frac{x}{1-y}),
			\end{equation*}
			for $x,y\in [0,1), x+y\in [0,1]$, with
			$f_\alpha(0) = f_\alpha(1) = 0$ and $f_\alpha(\frac{1}{2})=\frac{2^{\alpha-1}-1}{\alpha-1},\alpha \neq 1.$}
	\end{enumerate}
	The continuos analogue to Definition 3.1 is given as follows:
	\begin{definition}[Mathai's entropy (Continuous case)]
		\begin{equation}\label{eq1}
			M_\alpha(f) = \frac{\int_{-\infty}^{\infty} [f(x)]^{2-\alpha}dx-1}{\alpha-1}, \alpha \neq 1, -\infty < \alpha < 2,
		\end{equation}
		where $ f(x) $ is the density function.
	\end{definition}
For more properties and details see Mathai and Haubold (2006), 
 Sebastian (2015)
 .
\subsection{Mathai's Pathway Model }

By optimizing Mathai's entropy\index{entropy} measure, one can arrive at pathway model of Mathai (2005), which consists of  many of the standard
distributions in statistical literature as special cases. For fixed $\alpha$, consider the optimization of $M_\alpha(f)$, which implies optimization of $\int_x[f(x)]^{2-\alpha}{\rm d}x$, subject to the following conditions:
\begin{enumerate}[{(i)}]
\item $f(x)\geq 0,~\text{for all}~x$
\item $\int_x f(x){\rm d}x < \infty$
\item $\int_x x^{\rho(1-\alpha)}f(x){\rm d}x=\text{fixed for all}~f$
\item $\int_x x^{\rho(1-\alpha)+\delta}f(x){\rm d}x=\text{fixed for all}~f, \text{where}~\rho~ \text{and}~\delta ~\text{are fixed parameters}$
\end{enumerate}
By using calculus of variation, one can obtain the Euler equation as
\begin{eqnarray}
&&\frac{\partial}{\partial f}[f^{2-\alpha}-\lambda_1x^{\rho(1-\alpha)}f+\lambda_2x^{\rho(1-\alpha)+\delta}f]=0\nonumber\\
&&\Rightarrow (2-\alpha)f^{1-\alpha}=\lambda_1x^{\rho(1-\alpha)}[1-\frac{\lambda_2}{\lambda_1}x^\delta],\nonumber\\
&&\Rightarrow f_1=c_1x^\rho[1-a(1-\alpha)x^\delta]^{\frac{1}{1-\alpha}}
\end{eqnarray}
for $\frac{\lambda_2}{\lambda_1}=a(1-\alpha)$ for some $a>0$. For more details the reader may refer to the papers of Mathai and Haubold (2007, 2008).\\
\par
When $\alpha \rightarrow 1$,  the Mathai's entropy measure $M_{\alpha}(f)$ goes to
the Shannon \index{Shannon entropy} entropy\index{entropy} measure
and this is a variant of
 Havrda-Charv\'{a}t entropy\index{entropy}, and the variant form therein is Tsallis
 entropy\index{entropy}. Then when $\alpha$ increases from 1, $M_{\alpha }(f)$ moves away from Shannon entropy. Thus $\alpha$ creates a pathway moving from one function to another, through the generalized entropy also.

	For real scalar case the pathway model is the following:
	\begin{equation}\label{eq2.8}
		f_1(x) = cx^{\gamma-1}{[1-a(1-\alpha)x^{\delta}]}^{\frac{1}{1-\alpha}},
	\end{equation}
	$ a>0,\delta>0,1-a(1-\alpha)x^{\delta}>0,\gamma>0 $, where $ c $ is the normalizing constant and $ \alpha $ is the pathway parameter.
	For $ \alpha <1 $, (\ref{eq2.8}) remains as a generalized type-1 beta model. \\
	When $ \alpha >1 $, we can write $ 1- \alpha = -(\alpha-1), \alpha>1 $ so that (\ref{eq2.8}) is of the form
	\begin{equation}\label{eq2.9}
		f_2(x) = cx^{\gamma-1}{[1+a(\alpha-1)x^{\delta}]^{-\frac{1}{\alpha-1}}},x>0,
	\end{equation}
	which is a generalized type-2 beta model. When $ \alpha \rightarrow 1 $, (\ref{eq2.8}) and (\ref{eq2.9}) reduces to (using L'Hospital rule)
	\begin{equation}\label{eq2.10}
		f_3(x)=cx^{\gamma-1}e^{-ax^{\delta}}, x>0,
	\end{equation}
which is a generalized gamma model.\\
	When we analyse data of physical experiments and build models, we usually select a member from a parametric family of distributions. It is often found that the underlying distribution is in between two parametric family of distributions. In these situations, we can use Mathai's pathway model to create a pathway from one functional form to another.

\vskip.5cm{\section{Maximum entropy}}
\vskip.3cm
In this section we maximize Mathai's entropy under the conditions of a density function and fixed Generalized Gini index. A nonlinear ordinary first ordinary equation is obtained. The equation is further examined through three cases by applying conditions on the coefficients. The following theorem and remark are used to achieve this, for more details see Gelfand and Fomin (1963)
, Tanak et al. (2015)
.\\

\begin{thm}
		 Let $ L(y) $ be a functional of the form
		\begin{equation} \label{eq4.1}
			L(y) = \int_a^b G(y(x),y'(x),x)dx,
		\end{equation}
		where the given function G is continuous and has continuous first partial derivatives in each of its arguments. $ L(y) $ is defined on the set of functions $ y(x) $ which have continuous first derivatives in $ [a,b] $ and satisfy the boundary conditions $ y(a)=A, y(b)=B $. Then, a necessary condition for $ L(y) $ to have an extremum for a given function $ y(x) $, is that $ y(x) $ satisfy the Euler's equation:
		\begin{equation}\label{eq4.2}
			\frac{\partial G}{\partial y} - \frac{d}{dx}\frac{\partial G}{\partial y'} = 0.
		\end{equation}
	\end{thm}
	\textbf{Remark 4.1.} \textit{Suppose we are looking for an extremum in (\ref{eq4.1}) subject to the conditions $y(a) = A, y(b)=B$ and
		\begin{equation}\label{eq4.3}
			\int_{a}^{b} J_i(y(x),y'(x),x)dx = l_i, i=1, 2,\dots, m,
		\end{equation}
		where $l_1,l_2,\dots,l_m$ are constants. In this case, a necessary condition for an extremum is that
		\begin{equation}\label{eq4.4}
			\frac{\partial}{\partial y}(G+\sum_{i=1}^{m}\lambda_i J_i) - \frac{d}{dx}\frac{\partial}{\partial y'}(G+\sum_{i=1}^{m}\lambda_i J_i) = 0,
		\end{equation}
		where $\lambda_1,\lambda_2,\dots,\lambda_m$ are Lagrangian multipliers.}

	\subsection{Maximum Mathai's entropy under the constraints on generalized Gini Index}

	Let $X$ be a non-negative random variable with probability density function $f(x)$ and cummulative distribution function $F(x)$ such that
	\begin{equation}\label{eq4.5}
		\lim_{x\rightarrow\infty}f(x)=0.
	\end{equation}
	Here, we intend to find the distribution that maximizes the entropy (\ref{eq1}) subject to the constraints	
	\begin{enumerate}
		\item{$\int_{0}^{\infty}f(x)dx=1$}
		\item{$\int_{0}^{\infty}xf(x)dx=\mu$}
		\item{$G_\nu(F)=\delta$,}
	\end{enumerate}
	where
	\begin{equation}\label{eq6}
		G_\nu(F)=1-\frac{1}{\mu}\int_{0}^{\infty}{\bar{F}}^\nu(x)dx.
	\end{equation}
	Hence, the optimization problem is to maximize the functional
	\begin{equation}\label{eq7}
		M_\alpha(f) = \frac{1}{\alpha-1}\left[\int_{0}^{\infty}{[f(x)]}^{2-\alpha}dx-1\right], \,\, \alpha \neq 1, \, \alpha < 2,
	\end{equation}
	subject to the constraints
	\begin{enumerate}
		\item{$\int_{0}^{\infty}f(x)dx=1 $}
		\item{$\int_{0}^{\infty}xf(x)dx=\mu$}
		\item{$\int_{0}^{\infty}{\bar{F}}^\nu (x)dx=\eta$}.
	\end{enumerate}
	Now, (\ref{eq7}) can be rewritten as
	\begin{equation}\label{eq8}
		M_\alpha(f) = \frac{1}{\alpha-1}\left[\int_{0}^{\infty}[{[f(x)]}^{2-\alpha}-f(x)]dx\right].
	\end{equation}
	The Euler's equation is given by
	\begin{align}\label{eq9}
		&\frac{\partial}{\partial \bar{F}}\left[\frac{{[f(x)]}^{2-\alpha}-f(x)}{\alpha-1}-\lambda_1 f(x)-\lambda_2 xf(x)-\lambda_3{\bar{F}}^\nu (x)\right] \nonumber \\
		&+\frac{d}{dx}\frac{\partial}{\partial f}\left[\frac{{[f(x)]}^{2-\alpha}-f(x)}{\alpha-1}-\lambda_1f(x)-\lambda_2xf(x)-\lambda_3{\bar{F}}^\nu (x)\right] = 0.
	\end{align}
	where $\lambda_1, \lambda_2$ and $ \lambda_3 $ are the Lagrangian multipliers.
	From (\ref{eq9}), we have
	\begin{equation}\label{eq10}
		(2-\alpha){[f(x)]}^{1-\alpha}f'(x)+\lambda_2 f(x)+\lambda_3\nu{\bar{F}}^{\nu-1}(x)f(x)=0.
	\end{equation}
	\begin{equation}\label{eq11}
		\Rightarrow \frac{d}{dx}[{[f(x)]}^{2-\alpha}(x)]-\lambda_2\frac{d}{dx}\bar{F}(x)-\lambda_3 \frac{d}{dx}{\bar{F}}^\nu (x) = 0.
	\end{equation}
	Integrating (\ref{eq11}), we get,
	\begin{equation}\label{eq12}
		f^{2-\alpha}(x)-\lambda_2 \bar{F}(x)-\lambda_3{\bar{F}}^\nu (x)=c,
	\end{equation}
	where c is the constant of integration. By (\ref{eq4.5}) and $ \lim_{x \rightarrow \infty}\bar{F}(x)=0 $, equation(\ref{eq12}) becomes
	\begin{equation}\label{eq13}
		f^{2-\alpha}(x)-\lambda_2 \bar{F}(x)-\lambda_3{\bar{F}}^\nu (x)=0.
	\end{equation}
	\begin{equation}\label{eq14}
		\Rightarrow {\left[-\frac{d}{dx}\bar{F}(x)\right]}^{2-\alpha}-\lambda_2 \bar{F}(x)-\lambda_3{\bar{F}^{\nu}(x)}=0.
	\end{equation}
	Equation (\ref{eq14}) is a first order nonlinear ordinary differential equation. By applying conditions on the Lagrangian multipliers $ \lambda_2 $ and $ \lambda_3 $, the equation (\ref{eq14}) is examined through three cases.\\
	\hfill \\
	\textbf{Case 1:} When $\lambda_2 \neq0, \lambda_3 = 0$. \\
	Then (\ref{eq14}) becomes,
	\begin{align*}
		&{\left[-\frac{d}{dx}\bar{F}(x)\right]}^{2-\alpha}-\lambda_2 \bar{F}(x) &=& \, \,  0 \nonumber \\
		\Rightarrow & {\left[-\frac{d}{dx}\bar{F}(x)\right]}^{2-\alpha} &=& \, \,  \lambda_2 \bar{F}(x) \nonumber
	\end{align*}
	\begin{align}\label{eq14.1}
		\Rightarrow & -\frac{d}{dx}\bar{F}(x) &=& \, \,  {\left[\lambda_2 \bar{F}(x)\right]}^{\frac{1}{2-\alpha}} \nonumber \\
		\Rightarrow & \frac{d\bar{F}(x)}{\lambda_2^{\frac{1}{2-\alpha}}{[\bar{F}(x)]}^{\frac{1}{2-\alpha}}}&=& -dx \nonumber \\
		\Rightarrow & \frac{1}{\lambda_2^{\frac{1}{2-\alpha}}} {[\bar{F}(x)]} ^{\frac{-1}{2-\alpha}} d\bar{F}(x) &=& -dx \nonumber \\
		\Rightarrow & \frac{1}{\lambda_2 ^{\frac{1}{2-\alpha}}} \times \frac{{[\bar{F}(x)]}^{1-\frac{1}{2-\alpha}}}{1-\frac{1}{2-\alpha}} &=& -x + c, \nonumber \\
		\Rightarrow & \frac{2-\alpha}{\lambda_2^{\frac{1}{2-\alpha}}} \times \frac{{[\bar{F}(x)]}^\frac{2-\alpha - 1}{2-\alpha}}{2-\alpha - 1} &=& -x + c \nonumber \\
		\Rightarrow & \frac{2-\alpha}{\lambda_2^{\frac{1}{2-\alpha}}} \times \frac{{[\bar{F}(x)]}^{\frac{1-\alpha}{2-\alpha}}}{1-\alpha} &=& -x + c \nonumber \\
		\Rightarrow & \bar{F}(x) &=& \, \,  {\left[\frac{\lambda_2^{\frac{1}{2-\alpha}}(1-\alpha)(c-x)}{2-\alpha}\right]}^{\frac{2-\alpha}{1-\alpha}},
	\end{align}
	where c is the constant of integration. This case has been already dealt in Section 3.\\
	\hfill\\
	\textbf{Case 2:} When $\lambda_2 = 0, \lambda_3 \neq 0$. \\
	Then (\ref{eq14}) becomes,
	\begin{align*}\label{eq15}
		&{\left[-\frac{d}{dx}\bar{F}(x)\right]}^{2-\alpha}-\lambda_3 \bar{F}^\nu (x) &=& \, \,  0 \nonumber \\
		\Rightarrow & {\left[-\frac{d}{dx}\bar{F}(x)\right]}^{2-\alpha} &=& \, \,  \lambda_3 \bar{F}^\nu (x) \nonumber \\
		\Rightarrow & -\frac{d}{dx}\bar{F}(x) &=& \, \, {\left[\lambda_3 \bar{F}^\nu (x)\right]}^{\frac{1}{2-\alpha}} \nonumber \\
		\Rightarrow & \frac{d\bar{F}(x)}{\lambda_3^{\frac{1}{2-\alpha}}{[\bar{F}^\nu (x)]}^\frac{1}{2-\alpha}} &=& -dx \nonumber \\
		\Rightarrow & \frac{1}{\lambda_3^{\frac{1}{2-\alpha}}}{[\bar{F}(x)]}^{\frac{-\nu}{2-\alpha}}d\bar{F}(x) &=& -dx \nonumber \\
		\Rightarrow & \frac{2-\alpha}{\lambda_3^{\frac{1}{2-\alpha}}} \times \frac{{[\bar{F}(x)]}^\frac{2-\alpha - \nu}{2-\alpha}}{2-\alpha - \nu} &=& -x +c \nonumber
	\end{align*}
	\begin{align}
		\Rightarrow & \bar{F}(x) &=& \, \,  {\left[\frac{\lambda_3^{\frac{1}{2-\alpha}}(2-\alpha - \nu)(c-x)}{2-\alpha}\right]}^{\frac{2-\alpha}{2-\nu - \alpha}},
	\end{align}
	where c is the constant of integration.\\
	\hfill \\
	Suppose c = 0. Then (\ref{eq15}) becomes
	\begin{equation}\label{eq16}
		\bar{F}(x) = {\left[\frac{\lambda_3^{\frac{1}{2-\alpha}}(2-\alpha-\nu)(-x)}{2-\alpha}\right]^\frac{2-\alpha}{2-\nu-\alpha}}.
	\end{equation}
	At $ x = 0 $, $ \bar{F}(0) = 0 $ (using (\ref{eq16})) $ \neq 1 $. Therefore, c cannot be zero. \\
	\hfill \\
	Now, (\ref{eq15}) can be rewritten as
	\begin{equation} \label{eq17}
		\bar{F}(x) = {\left[\frac{\lambda_3^{\frac{1}{2-\alpha}}(2-\alpha-\nu)c(1-\frac{x}{c})}{2-\alpha}\right]}^\frac{2-\alpha}{2-\nu-\alpha}.
	\end{equation}
	Using $ \bar{F}(0) = 1 $, we get from (\ref{eq17}),
	\begin{equation} \label{eq18}
		{\left[\frac{\lambda_3^{\frac{1}{2-\alpha}}(2-\alpha - \nu)c}{2-\alpha}\right]}^{\frac{2-\alpha}{2-\nu-\alpha}} = 1.
	\end{equation}
	Therefore, (\ref{eq17}) reduces to
	\begin{equation} \label{eq19}
		\bar{F}(x) = {\left(1-\frac{x}{c}\right)}^{\frac{2-\alpha}{2-\nu-\alpha}}.
	\end{equation}
	Now using property $ \lim\limits_{x \to \infty} \bar{F}(x) = 0 $, $ \lim\limits_{x \to \infty} {\left(1-\frac{x}{c}\right)}^{\frac{2-\alpha}{2-\nu-\alpha}}  = 0 $ (see \ref{eq18}), only if,
	\begin{equation} \label{eq20}
		2< \nu + \alpha.
	\end{equation}
	Using property, $ f(x) = -\frac{d}{dx}\bar{F}(x) $, where $ f(x) $ is the density function, differentiating (\ref{eq19}), we get
	\begin{equation} \label{eq21}
		f(x) = {\left(\frac{1}{c}\right)}{\left(\frac{2-\alpha}{2-\nu-\alpha}\right)}{\left(1-\frac{x}{c}\right)}^{\frac{\nu}{2-\nu-\alpha}}.
	\end{equation}
	Using property, $ \int_{0}^{\infty} f(x) = 1 $, from (\ref{eq21}), we get $ c = \pm 1 $, . Since $ f(x) \geq 0 $, we get,  c = -1.\\
	Thus, the density function is given by,
	\begin{equation} \label{eq22}
		f(x) = {\left(\frac{\alpha-2}{2-\nu-\alpha}\right)}{\left(1+x\right)}^{\frac{\nu}{2-\nu-\alpha}},
	\end{equation}
	where $ x \geq 0, -\infty< \alpha <2, \alpha \neq 1, \nu > 1 $ and $ \nu + \alpha >2 $.\\
	The density function of Lomax distribution is given by,
	\begin{equation} \label{eq23}
		f(x) =\frac{\beta}{\lambda}{\left(1+\frac{x}{\lambda}\right)}^{-(\beta+1)},  x \geq 0, \beta > 0 \, \text{and} \, \lambda > 0.
	\end{equation}
	Now, (\ref{eq22}) can be rewritten as
	\begin{equation} \label{eq24}
		f(x) = {\left(\frac{\alpha-2}{2-\nu-\alpha}\right)}{(1+x)}^{-(\frac{-\nu}{2-\nu-\alpha})}.
	\end{equation}
	Comparing, (\ref{eq23}) and (\ref{eq24}), we get $ \beta = \frac{\alpha-2}{2-\nu-\alpha} $. Therefore, (\ref{eq24}) becomes,
	\begin{equation} \label{eq25}
		f(x) = \beta {(1+x)}^{-(\beta +1)}.
	\end{equation}
	Thus (\ref{eq22}) is the density of Lomax distribution with parameters $ \lambda=1 $ and $ \beta = \frac{\alpha-2}{2-\nu-\alpha} $. \\
	\hfill \\
	\textbf{Case 3:} When $\lambda_2 \neq 0, \lambda_3 \neq 0$. \\
	In this case, we do not have an explicit solution for (\ref{eq14}).

\section{Maximization of Mathai's entropy under the constriants on Gini mean difference index}
Here, we intend to find the distribution that maximizes the entropy (\ref{eq1}) subject to the constraints
\begin{enumerate}
	\item $ \int_{0}^{\infty} f(x) dx = 1 $
	\item $ \int_{0}^{\infty} xf(x) dx = \mu $
	\item $ D(f) = \gamma  $.
\end{enumerate}
Using (\ref{gp2}), the optimization problem is to maximize the functional
\begin{equation}\label{gp3}	
M_{\alpha} (f) = \frac{1}{\alpha -1} \left[\int_{0}^{\infty} {[f(x)]}^{2-\alpha} dx -1 \right], \, \alpha \neq 1, \, \alpha <2, \\
\end{equation}
subject to the constraints
\begin{enumerate}
	\item $ \int_{0}^{\infty} f(x) dx = 1 $
	\item $ \int_{0}^{\infty} xf(x) dx = \mu $
	\item $ \int_{0}^{\infty} {\bar{F}}^2 (x) dx = \phi  $.
\end{enumerate}
If we proceed as in the previous section, we get (24) with $ \nu = 2 $. The resulting non-linear differential equation is given by,
\begin{equation}\label{gp4}
{\left[-\frac{d}{dx}\bar{F}(x)\right]}^{2 - \alpha} - \lambda_2 \bar{F}(x) - \lambda_3 \bar{F}^2(x)=0.
\end{equation}
\hfill \\
\textbf{Case 1:} When $ \lambda_2 \neq 0, \lambda_3 =0 $.\\
In this case, we obtain the same expression as given in (\ref{eq14.1}). \\
\hfill \\
\textbf{Case 2:} When $ \lambda_2 = 0, \lambda_3 \neq 0 $.\\
Proceeding as in the previous section, the density function is obtained by substituting $ \nu =2 $ in (\ref{eq22}). The resulting density function is,
\begin{equation}\label{gp5}
f(x) = \left(\frac{2-\alpha}{\alpha}\right){(1+x)}^{-\frac{2}{\alpha}},
\end{equation}
where $ x \geq 0, \alpha>0 $ and $ \alpha \neq 1 $. Equation (\ref{gp5}) is the special case of Lomax distribution given in equation (\ref{eq23}) with parameters $ \lambda = 1 $ and $ \beta = \frac{2-\alpha}{\alpha}$.\\
\hfill \\
\textbf{Case 3:} When $ \lambda_2 \neq 0, \lambda_3 \neq 0 $.\\
In this case, we do not have an explicit solution for (\ref{gp4}).

\vskip.5cm{\section{Numerical Illustration}}
\vskip.3cm

An application of the distribution in (\ref{eq22}) is presented here by comparing the fits of this model to some selected one-parameter distributions available in the literature.\\
	The density function and distribution function of exponential distribution is
	\begin{gather*} \label{5.1}
f(x) = \theta e^{-\theta x}\\
F(x) = 1-e^{-\theta x}
	\end{gather*}
where $ x\geq 0, \, \theta >0.  $ \\
The density function and distribution function of Lindley distribution (Lindley (1958), Hafez et. al. (2020)) is
\begin{gather*} \label{5.12}
f(x) = \frac{\theta^2}{1+\theta}(1+x) e^{-\theta x} \\
F(x) = 1 - \frac{(1+\theta+\theta x)e^{-\theta x}}{1+\theta}
	\end{gather*}
where $ x\geq 0, \, \theta >0.  $ \\
The density function and distribution function of Akash distribution (Shankar (2015)) is
\begin{gather*} \label{5.13}
	f(x) = \frac{\theta ^3}{\theta ^2 +2} (1+x^2) e^{-\theta x}\\
	F(x) = 1-\left[1+\frac{\theta x (\theta x +2)}{\theta^2 +2}\right]e^{-\theta x}
	\end{gather*}
where $\, \,\, \, x\geq 0, \, \theta >0.  $ \\
The density function and distribution function of Pranav distribution (Shukla (2018)) is
\begin{gather*}\label{5.14}
	f(x) = \frac{\theta ^4}{\theta ^4 +6} (\theta + x^3) e^{-\theta x}\\
	F(x) = 1-\left[1+\frac{\theta x (\theta^2 x^2 + 3\theta x +6)}{\theta^4 +6}\right]e^{-\theta x}
	\end{gather*}
where $ x\geq 0, \, \theta >0.  $ \\
The density function and distribution function of Ishitha distribution (Shankar and Shukla (2017)) is
\begin{gather*} \label{5.15}
	f(x) = \frac{\theta ^3}{\theta ^3 +2} (\theta + x^2) e^{-\theta x}\\
	F(x) = 1-\left[1+\frac{\theta x (\theta x + 2)}{\theta^3 +2}\right]e^{-\theta x}
	\end{gather*}
where $ x\geq 0, \, \theta >0.  $ \\
The density function and distribution function of Ram Awadh distribution (Shukla (2018)) is
\begin{gather*} \label{5.16}
	f(x) = \frac{\lambda ^6}{\lambda ^6 +120} (\lambda + x^5)e^{-\lambda x}\\
	F(x) = 1-\left[1+\frac{\lambda x (\lambda^4x^4+5\lambda^3x^3+20\lambda^2x^2+60\lambda x+120)}{\lambda^6 +120}\right]e^{-\lambda x}
	\end{gather*}
where $ x\geq 0, \, \lambda >0.  $ \\
The density function and distribution function of Sujatha distribution (Shankar (2016)) is
\begin{gather*} \label{5.17}
	f(x) = \frac{\theta ^3}{\theta ^2+ \theta +2} (1+x+x^2) e^{-\theta x}\\
	F(x) = 1-\left[1+\frac{\theta x (\theta x + \theta + 2)}{\theta^2 + \theta +2}\right]e^{-\theta x}
	\end{gather*}
where $ x\geq 0, \, \theta >0.  $ \\
	Here, the distribution (\ref{eq22}) with Gini index parameter $ \nu=3 $ and  various one-parameter distributions mentioned above is fitted to the time series data of loss ratios (yearly data) for eathquake insurance in California from 1971 till 1994 ((Jaffee and Russell) (1997)). The data are: 17.4, 0.0, 0.6, 3.4, 0.0, 0.0, 0.7, 1.5, 2.2, 9.2, 0.9, 0.0, 2.9, 5.0, 1.3, 9.3, 22.8, 11.5, 129.8, 47.0, 17.2, 12.8, 3.2, 2272.7\\
	\hfill \\
Table 5.1: Maximum Likelihood Estimates (MLE), Kolmogorov Smirnov (KS) (and its p value), Akaike Information Criterion (AIC) and Bayesian Information Criterion (BIC)
\begin{center}
\begin{tabular}{|c|c|c|c|c|c|c|}
\hline
Model & MLE & KS & P value & AIC & BIC\\
\hline
\textbf{Distribution } & \textbf{0.9697608} & \textbf{0.16667} &\textbf{ 0.5176} & \textbf{172.8832} & \textbf{174.0612} \\

\textbf{considered} &  &  &  &  &   \\
Exponential & 0.009333437 & 0.68343 & 3.666 $ \times 10^{-10} $ & 274.3593 & 275.5373 \\
Lindley & 0.01849737 & 0.8025 & 7.516 $ \times 10^{-14} $ & 389.2868 &390.4649 \\
Akash & 0.027993 & 0.84776 & 2.109 $ \times 10^{-15} $ & 523.8548 &525.0329 \\
Pranav & 0.0390226 & 0.86204 & 6.661 $ \times 10^{-16} $ & 695.1357 & 696.3137 \\
Ishitha & 0.02963961 & 0.84373 & 2.887 $ \times 10^{-15} $ & 561.0583 & 562.2363 \\
Ram Awadh & 0.05786601 & 0.8726 & 2.22 $ \times 10^{-16} $ & 968.356 & 969.534 \\
Sujatha & 0.02786496 & 0.84661 & 2.331 $ \times 10^{-15} $ & 517.0664 & 518.2444\\
\hline
\end{tabular}
\end{center}
Thus the distribution (\ref{eq22}) provides a better fit to the given data.\\
\hfill \\

\vskip.5cm{\section{Conclusion}}
\vskip.3cm
Mathai's entropy, one of the several generalizations of Shannon entropy, can be maximized under different sets of constraints to obtain density functions. Here, Mathai's entropy is maximized, subject to the constraints of a given generalized Gini index and the conditions of a density function. The resulting distribution is a special case of Lomax distribution.  Mathai's entropy maximization based on Gini mean
difference constraints showed that the maximum Mathai's entropy distributions is the special case of Lomax distribution
discussed in previous case.

\end{document}